\documentclass[preprint,prl,nofootinbib]{revtex4-1}
\usepackage{graphicx}
\usepackage{dcolumn}
\usepackage{amsmath,amsthm,amssymb}

\begin{document}
\title{Shrinkers, expanders, and the unique continuation beyond generic blowup in the heat flow for harmonic maps between spheres}

\author{Pawe\l {}  Biernat}
\affiliation{M. Smoluchowski Institute of Physics, Jagiellonian University, Krak\'ow, Poland}
\author{Piotr Bizo\'n}
\affiliation{M. Smoluchowski Institute of Physics, Jagiellonian University, Krak\'ow, Poland}

\date{\today}
\begin{abstract}
Using mixed analytical and numerical methods we investigate the development of singularities in the heat flow for corotational harmonic maps from the $d$-dimensional sphere to itself for $3\leq d\leq 6$.  By gluing together shrinking and expanding asymptotically self-similar solutions we construct global weak solutions which are smooth everywhere except for a sequence of times $T_1<T_2<\dots <T_k<\infty$ at which there occurs the type I blow-up at one of the poles of the sphere.
We give evidence that in the generic case the continuation beyond blow-up is unique, the topological degree of the map changes by one at each blow-up time $T_i$, and eventually the solution comes to rest at the zero energy constant map.
\end{abstract}

\maketitle
\addtolength{\topmargin}{-0.8pc}
\addtolength{\textheight}{1.6pc}
\section{Introduction}
\noindent
Let $M$ and $N$ be Riemannian manifolds with metric tensors $g_{ij}$ and $G_{AB}$ in some local coordinates $\{x^i\}$ and  $\{X^A\}$. A map $X:M\rightarrow N$ is called harmonic if
it is a critical point of the energy
\begin{equation}\label{energy}
    E(X)= \int_M e(X) \sqrt{g}\, dx\,,\qquad
    e(X)=\frac{1}{2} \frac{\partial X^A}{\partial x^i} \frac{\partial X^B}{\partial x^j} G_{AB}\, g^{ij}\,.
\end{equation}
 In this paper we consider harmonic maps from the $d$-dimensional unit sphere to itself, i.e. $M=N=S^d$ with $g_{ij}$ and $G_{AB}$ being standard round metrics. We parametrize $S^d$ by spherical coordinates $(\theta,\phi)$, where $\theta$ is colatitude ($0\leq \theta\leq \pi$) and $\phi$ is a point on the equator $S^{d-1}$ of $S^d$. We restrict our attention to corotational maps of the form $(\theta,\phi)\rightarrow (U(\theta),\phi)$. For such maps we have
\begin{equation}\label{E(U)}
    E(U)=\frac{1}{2} \int_0^{\pi} \left(U_{\theta}^2+(d-1) \, \frac{\sin^2{\!U}}{\sin^2{\!\theta}}\right)
\sin^{d-1}{\!\theta}\, d\theta\,,
\end{equation}
where for convenience we dropped the multiplicative factor $vol(S^{d-1})$ coming from the integration over $\phi$. The Euler-Lagrange equation corresponding to the energy \eqref{E(U)} reads
\begin{equation}\label{eq(u)}
    \frac{1}{\sin^{d-1}{\!\theta}} \left(\sin^{d-1}{\!\theta}\, U_{\theta} \right)_{\theta}- \frac{d-1}{2} \frac{\sin(2U)}{\sin^2{\!\theta}}=0\,.
\end{equation}
It was shown in \cite{bc} that for $3\leq d\leq 6$ Eq.\eqref{eq(u)} has a countable sequence $\{U_n\}$ of smooth solutions of degree zero and one. These solutions may be viewed as excitations of the ground states: the constant map $U_0=0$ and the identity map $U_1=\theta$, for even and odd values of $n$, respectively (for $d\geq 7$ these excitations disappear). For $n\rightarrow \infty$ the solutions $U_n(\theta)$ converge (nonuniformly)
to the (singular) equator map $U_{\infty}=\pi/2$.
Later,  Corlette and Wald \cite{cw} rederived and extended these results using Morse theory methods. Their approach helped to identify the two key features which are responsible for the existence of infinitely many solutions: the presence of the antipodal reflection symmetry $U\rightarrow \pi-U$ and the existence of the singular map $U_{\infty}=\pi/2$ of infinite index which is invariant under this symmetry.  An essential ingredient of the Morse theoretic argument is an energy decreasing flow in the space of maps. In \cite{cw} this flow was defined in a somewhat \emph{ad hoc} manner to ensure that it has all the desired technical properties. One might wonder if it is possible to repeat the Corlette-Wald argument using the ordinary heat flow. This would be interesting, for instance, in numerical implementations of the argument for similar systems. The main technical difficulty is that the heat flow can develop singularities in finite time.
If this happens, in order to save the argument,
  one must find a way to continue the flow past a singularity in a unique manner.
Although an analysis of this issue was the original motivation for this paper, the problem of uniqueness of continuation beyond blow-up in the heat flow for harmonic maps seems interesting in its own right, regardless of possible applications to elliptic problems.

 The aim of this paper is two-fold. First, we describe the precise asymptotics of blow-up in the heat flow for corotational harmonic maps from the $d$-dimensional sphere to itself for $3\leq d\leq 6$. We show that blow-up has the form of a shrinking self-similar solution (shrinker, for short). In turns out that among infinitely many shrinkers (whose existence was proved by Fan \cite{f}), there is exactly one which is linearly stable.
  We provide numerical evidence that this stable shrinker  determines the generic profile of blow-up. Second, we continue the flow past the singularity by gluing a suitable expanding self-similar solution (expander, for short). We find that there is exactly one expander which can be glued to the stable shrinker and consequently the continuation beyond  the generic blow-up  is unique.

 The scenario of incomplete blow-up and self-similar global "peaking solutions", that is solutions which shrink self-similarly, blow up, and then expand self-similarly for a while (with this scenario possibly repeating a number of times) has been studied in the past  for the harmonic map flow \cite{i} and other parabolic equations: the semilinear heat equations \cite{lt, gv}, the mean curvature flow \cite{h,aci,aiv,i}, the Yang-Mills flow \cite{gas},
 the Ricci flow \cite{fik}, and more recently for a fourth-order reaction-diffusion equation \cite{gal2}. Most of these studies emphasized non-uniqueness of continuation beyond blow-up. To our knowledge, this is the first work which demonstrates (by heuristic and numerical means) that for the generic blow-up the continuation is unique. As we shall see below, the uniqueness of continuation is contingent upon certain \emph{quantitative}
 properties of self-similar solutions and thus may be hard to prove.

 The rest of the paper is organized as follows. In section~2 we introduce the  heat flow for equivariant harmonic maps from $S^d$ (or $\mathbb{R}^d$) into $S^d$ and recall basics facts about blow-up. Section~3 is devoted to self-similar solutions of the heat flow for harmonic maps from $\mathbb{R}^d$ to $S^d$. Using matched asymptotics we derive asymptotic scaling formulae for the parameters of self-similar solutions. In section~4 we analyze the linear stability of self-similar solutions. In section~5 we study the continuation beyond blow-up and formulate the main result of this paper, that is the conjecture about the uniqueness of continuation in the generic case. Numerical evidence supporting this conjecture is presented in section~6.
  Finally, in section~7 we indicate possible extensions of our results.
 \addtolength{\topmargin}{+0.pc}
\addtolength{\textheight}{-1.6pc}
\section{Preliminaries}
We consider the heat flow equation
\begin{equation}\label{heat}
  U_t=  \frac{1}{\sin^{d-1}{\!\theta}} \left(\sin^{d-1}{\!\theta}\, U_{\theta} \right)_{\theta}- \frac{d-1}{2} \frac{\sin(2U)}{\sin^2{\!\theta}}\,,
\end{equation}
with initial and boundary conditions
\begin{eqnarray}\label{ic}
    U(0,\theta)&=&h(\theta)\in C^{\infty}[0,\pi]\,,\\
    U(t,0) &=&h(0)=0\,,\\
    U(t,\pi)&=&h(\pi)=k\pi\,,
\end{eqnarray}
where an integer $k$ is the topological degree of the map.  As long as the flow is smooth, the solution remains in the given homotopy class (i.e., the degree $k$ does not change). It follows from \eqref{heat} that for a smooth solution there holds
\begin{equation}\label{dEdt}
    \frac{dE}{dt}=-\int_0^{\pi} U_t^2 \sin^{d-1}{\!\theta} \,d\theta\,,
\end{equation}
which shows that Eq.\eqref{heat} is the gradient flow for the energy \eqref{E(U)}.
Thus, one might  expect that for $t\rightarrow\infty$ the solution $U(t,\theta)$ will converge to a critical point of $E$, i.e. a harmonic map. Unfortunately, as mentioned in the introduction, this expectation is too naive because in general the flow develops  singularities in finite time.
Indeed, it follows from general results for harmonic maps between  compact manifolds  that for any initial map with nonzero degree and sufficiently small energy the solution must blow up in finite time (see Thm 1.12 in \cite{s1}).

 For equation \eqref{heat}, by symmetry, the singularity must occur at one of the poles. Since the blow-up is a localized phenomenon, the curvature of the domain manifold plays no role in the description of asymptotics of blow-up.  Thus, from here until section~6  we replace the domain $S^d$ by its tangent space at the pole, $\mathbb{R}^d$, and consider the heat flow for corotational harmonic maps from $\mathbb{R}^d$ to $S^d$
  \begin{equation}\label{heat_flat}
  u_t=  \frac{1}{r^{d-1}} \left(r^{d-1} u_r \right)_r- \frac{d-1}{2 r^2} \sin(2u)\,,
\end{equation}
   where $u=u(t,r)$ ($r=|x|$).  Such maps enjoy scale invariance: if $u(t,r)$ is a solution, so is $u_{\lambda}(t,r)=u(t/\lambda^2,r/\lambda)$ for any positive number $\lambda$. Solutions which are invariant under rescaling, that is $u_{\lambda}=u$, are called self-similar. The self-similar solutions play the key role in the dynamics of  type I blow-up\footnote{It is customary to divide singularities into two types: a singularity for which $(T-t)|\nabla u|^2$ is bounded as $t\nearrow T$ is said to be of type I; otherwise it is said to be of type II.} so the next three sections are devoted to their existence and properties.

   Throughout the rest of this paper we assume that $3\leq d\leq 6$.
\section{Self-similar solutions}
\subsection{Shrinkers}
\noindent
Let us assume that a solution of Eq.\eqref{heat_flat} develops a type I singularity at $r=0$ in a finite time $T$, i.e. $(T-t) u^2_r(t,0)$ is bounded as $t\nearrow T$.
To describe the formation of the singularity it is convenient to introduce new variables
\begin{equation}
s=-\ln(T-t),\,\,\,y=\frac{r}{\sqrt{T-t}}\,,\quad f(s,y)=u(t,r)\,.
\end{equation}
In these variables Eq.\eqref{heat_flat} takes the form
\begin{equation}\label{eq_f}
  f_s =
  \frac{1}{\rho} \left(\rho
  f_y\right)_y -\frac{d-1}{2 y^2}\,\sin(2f)\,,\qquad \rho(y)=y^{d-1} \exp(-y^2/4)\,.
  \end{equation}
This equation can be viewed as the gradient flow for the  functional
\begin{equation}\label{e}
    \mathcal{E}(f)=\frac{1}{2} \int_0^{\infty}  \left(f_y^2+\frac{d-1}{y^2} \sin^2\!{f} \right)\, \rho\, dy\,.
\end{equation}
We shall refer to $\mathcal{E}(f)$ as the conformal energy because it is the energy for maps from $(\mathbb{R}^d,e^{-\frac{y^2}{2(d-2)}} \delta$) to $S^d$.
A simple calculation gives
\begin{equation}\label{de}
    \frac{d\mathcal{E}}{ds}= - \int_0^{\infty} f_s^2\, \rho \,dy\,,
\end{equation}
hence the conformal  energy is monotonically decreasing. The assumption that the blow-up is of type I implies that  $f_y$  is uniformly bounded as $s\rightarrow \infty$, hence the flow must converge to a critical point of the conformal energy, that is a solution of the Euler-Lagrange equation $\delta \mathcal{E}(f)=0$ for $f(y)$
\begin{equation}\label{shrink_eq}
    f''+\left(\frac{d-1}{y}-\frac{y}{2}\right)\,f'-\frac{d-1}{2 y^2} \sin(2f)=0\,.
\end{equation}
  Note that an endpoint of evolution cannot be the trivial solution $f=0$ as this would contradict the occurrence of blow-up at time $T$. Thus, the study of type I blow-up reduces to the study of nonconstant solutions of Eq.\eqref{shrink_eq}.
We shall call such solutions shrinkers.

Let us discuss now existence and properties of shrinkers.
Regular solutions of Eq.\eqref{shrink_eq} behave near $y=0$ as follows
\begin{equation}\label{a}
    f(y)=a y -\frac{a(4 d a^2-4 a^2-3)}{12 (2+d)}\, y^3 +\mathcal{O}(y^5)\,,
\end{equation}
where $a$ is a free parameter. Regular solutions at infinity behave as
\begin{equation}\label{b}
    f(y)=\frac{\pi}{2}+ b-\frac{(d-1) \sin(2b)}{2 y^2}+\mathcal{O}(y^{-4})\,,
\end{equation}
where $b$ is a free parameter. Using a shooting method Fan \cite{f} proved that for $3\leq d\leq 6$ there is an infinite sequence of  pairs $(a_n,b_n)$ for which the local solutions \eqref{a} and \eqref{b}  are smoothly connected by a globally regular solution $f_n(y)$. The integer index $n$ denotes the number of intersections of the solution $f_n(y)$ with $\pi/2$ (see Figure~1). As $n\rightarrow \infty$ the shrinkers converge (nonuniformly)
to the equator map $f_{\infty}=\pi/2$ and correspondingly $E(f_n)\rightarrow E(f_{\infty})=2^{d-1}\Gamma(\frac{d-1}{2})$. Some quantitative characteristics of shrinkers are displayed in Table~I.
\begin{table}[h]
  \begin{tabular}{|c|c|c|c|}\hline
$n$      & $a_n$                       & $b_n$                         & $E_n$          \\\hline
$\,1\,$  & $\,2.738753\,$              & $\,0.573141\,$                & $\,1.485688\,$ \\\hline
$\,2\,$  & $\,2.927644 \cdot 10^{1}\,$ & $\,-0.184519\,$               & $\,1.738165\,$ \\\hline
$\,3\,$  & $\,3.141830 \cdot 10^{2}\,$ & $\,0.566142 \cdot 10^{-1}\,$  & $\,1.771588\,$ \\\hline
$\,4\,$  & $\,3.376630 \cdot 10^{3}\,$ & $\,-0.172776 \cdot 10^{-1}\,$ & $\,1.776470\,$ \\\hline
$\,5\,$  & $\,3.629513 \cdot 10^{4}\,$ & $\,0.527011 \cdot 10^{-2}\,$  & $\,1.778116\,$ \\\hline
$\,6\,$  & $\,3.901390 \cdot 10^{5}\,$ & $\,-0.160744 \cdot 10^{-2}\,$ & $\,1.779706\,$ \\\hline
$\,7\,$  & $\,4.193637 \cdot 10^{6}\,$ & $\,0.490287 \cdot 10^{-3}\,$  & $\,1.781650\,$ \\\hline
$\,8\,$  & $\,4.507777 \cdot 10^{7}\,$ & $\,-0.149542 \cdot 10^{-3}\,$ & $\,1.784095\,$ \\\hline
$\,9\,$  & $\,4.845449 \cdot 10^{8}\,$ & $\,0.456120 \cdot 10^{-4}\,$  & $\,1.787199\,$ \\\hline
$\,10\,$ & $\,5.208415 \cdot 10^{9}\,$ & $\,-0.139121 \cdot 10^{-4}\,$ & $\,1.791128\,$ \\\hline
\end{tabular}


  \caption{Parameters of the first ten shrinkers for $d=3$.
    \label{table:energy}}
\end{table}
\begin{figure}[ht]
\center
\includegraphics[width=0.7\textwidth]{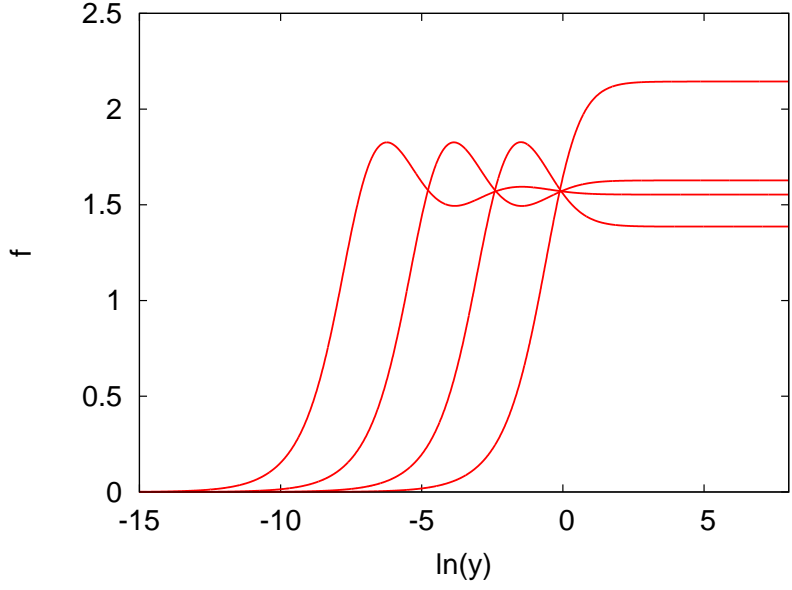}
\caption{The profiles of the first four shrinkers for $d=3$.
\label{fig:shrinkers}}
\end{figure}
\vskip 1cm
From the
 shooting argument in \cite{f} it follows that $a_n\rightarrow \infty$ and $b_n\rightarrow 0$ as $n\rightarrow\infty$.
 We shall now use this fact to describe the behaviour of shrinkers for large $n$ (we follow here a similar argument given in \cite{bc}).
 Let $\xi=a y$ and $\phi(\xi)=f(y)$. In terms of these variables Eq.\eqref{shrink_eq} becomes
\begin{equation}\label{phi}
   \phi''+ \left(\frac{d-1}{\xi}-\frac{\xi}{2 a^2}\right)\,\phi' -\frac{d-1}{2 \xi^2} \sin(2\phi)=0
\end{equation}
with the initial condition  $\phi(\xi)\sim \xi$ near $\xi=0$. For $a\rightarrow \infty$, solutions of this equation tend uniformly on any compact interval to solutions of the limiting equation
\begin{equation}\label{philim}
    \tilde \phi''+ \frac{d-1}{\xi}\,\tilde \phi' -\frac{d-1}{2 \xi^2} \sin(2\tilde \phi)=0
\end{equation}
with the same initial condition $\tilde \phi(\xi)\sim \xi$ near $\xi=0$.  Using the standard phase-plane analysis we get for $1\ll\xi\ll a$
\begin{equation}\label{asym}
    \tilde \phi(\xi)\simeq \frac{\pi}{2} + \alpha\, \xi^{-\frac{d-2}{2}} \sin(\omega \ln{\xi}+\delta)\,,\qquad \omega=\frac{\sqrt{8d-d^2-8}}{2}\,,
\end{equation}
where the amplitude $\alpha$ and the phase $\delta$ are uniquely determined by the initial condition $\tilde\phi'(0)=1$.
Returning to the original variables we obtain for $1/a\ll y\ll 1$
\begin{equation}\label{asymf}
    f(y)\simeq\frac{\pi}{2} + a^{-\frac{d-2}{2}} \alpha\, y^{-\frac{d-2}{2}} \sin(\omega \ln{y} +\omega \ln{a}+\delta)\,.
\end{equation}
On the other hand, for $y\gg1/a$ the solution is close to $\pi/2$ so we can write
\begin{equation}\label{scal2b}
    f(y)\simeq \frac{\pi}{2} + b\, h(y)\,,
\end{equation}
where $h(y)$ is the solution of the linearized equation
\begin{equation}\label{h_eq}
    h''+\left(\frac{d-1}{y}-\frac{y}{2}\right) h' +\frac{d-1}{y^2}\, h=0
\end{equation}
normalized by the condition $h(\infty)=1$. For $1/a\ll y\ll 1$ we have
\begin{equation}\label{h}
h(y)\simeq \alpha_1 y^{-\frac{d-2}{2}} \sin(\omega \ln{y}+\delta_1)\,,
\end{equation}
where $\alpha_1$ and $\delta_1$ are uniquely determined by the normalization condition
$h(\infty)=1$.
Using \eqref{h} and matching  the solutions (\ref{asymf}) and (\ref{scal2b}) we get
\begin{equation}\label{match}
  a^{-\frac{d-2}{2}} \alpha \sin(\omega \ln{y}+\omega \ln{a} + \delta)\simeq b\, \alpha_1 \sin(\omega \ln{y}+\delta_1)\,,
\end{equation}
hence
\begin{equation}
\omega \ln{a}+\delta\simeq \delta_1+n\pi,\qquad b\simeq (-1)^n \frac{\alpha}{\alpha_1} a^{-\frac{d-2}{2}},
\end{equation}
 which yields the scaling laws for large $n$
\begin{equation}\label{scala}
    a_n \simeq C \exp{\left(\dfrac{n\pi}{\omega}\right)},\qquad b_n \simeq (-1)^n D \exp{\left(-\dfrac{n (d-2)\pi}{2\omega}\right)},
\end{equation}
where $C=\exp[(\delta_1-\delta)/\omega]$ and $D=\dfrac{\alpha}{\alpha_1} C^{-\frac{d-2}{2}}$. Numerical parameters of shrinkers are displayed in Table~1 and compared with the asymptotic expressions \eqref{scala} in Fig.~2.

\begin{figure}[h]
\center
\includegraphics[width=0.48\textwidth]{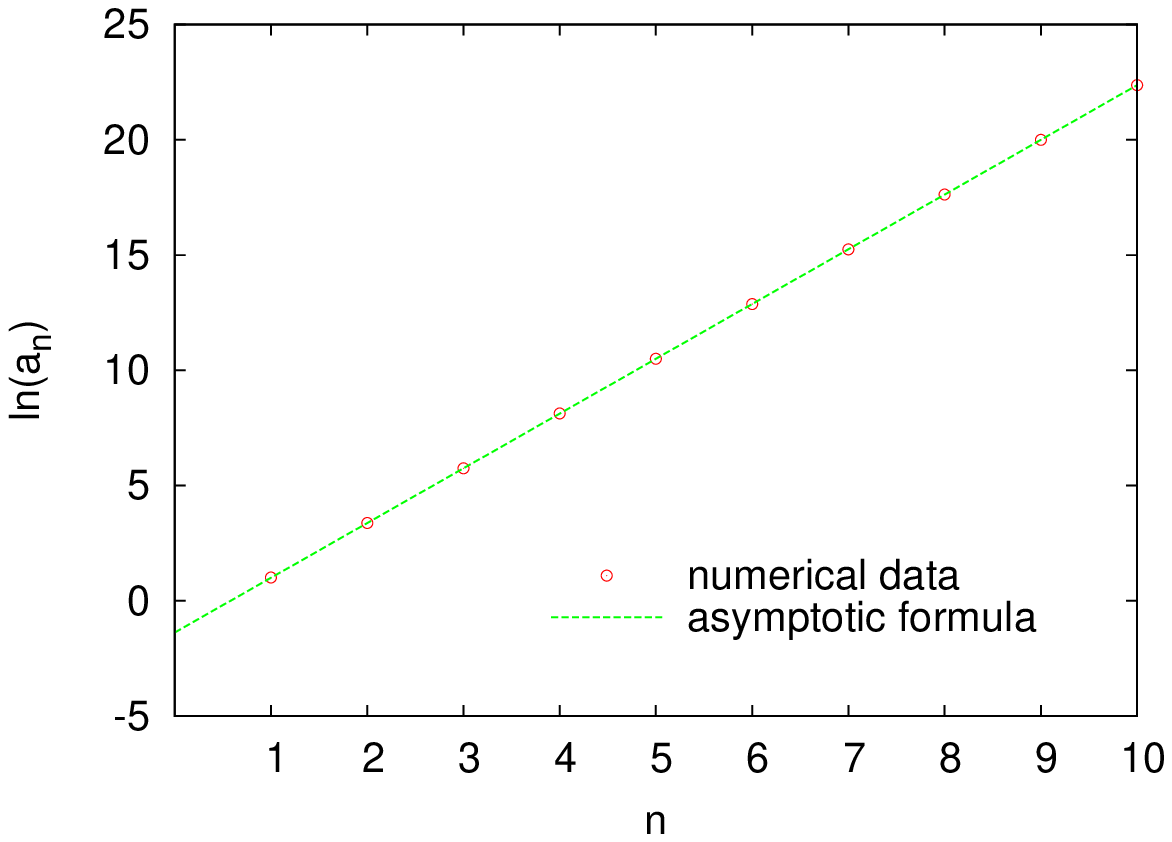}
\includegraphics[width=0.48\textwidth]{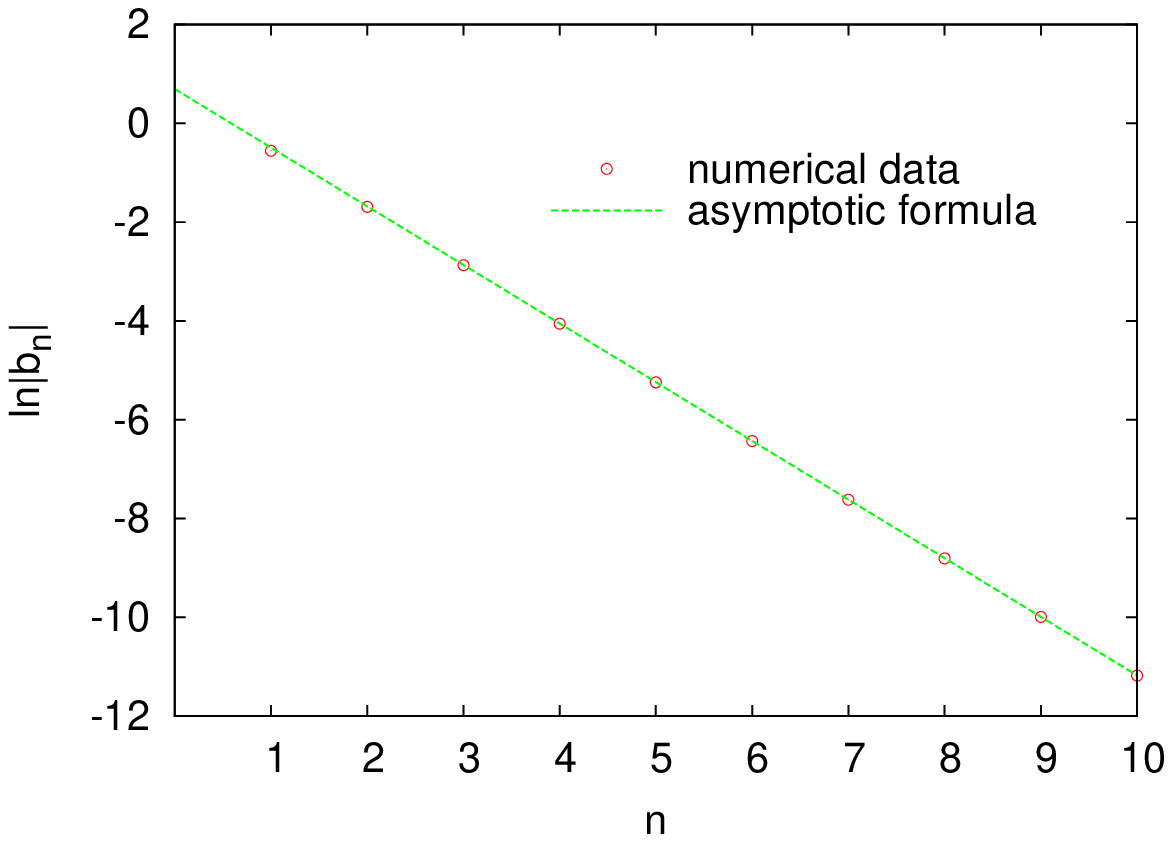}
\caption{The asymptotic formulae \eqref{scala} are shown to give excellent approximations for the parameters of shrinkers even for small $n$ (here $d=3$).
\label{fig:scaling}}
\end{figure}
\subsection{Expanders}
\noindent To describe the behaviour of solutions for $t>T$ we  introduce new variables
\begin{equation}
\sigma=\ln(t-T),\,\,\,y=\frac{r}{\sqrt{t-T}}\,,\quad F(\sigma,y)=u(t,r)\,,
\end{equation}
in which Eq.\eqref{heat_flat} takes the form
\begin{equation}\label{eq_F}
  F_{\sigma} =
  \frac{1}{R} \left(R\,
  F_y\right)_y -\frac{d-1}{2 y^2}\,\sin(2F)\,,\qquad R(y)=y^{d-1} \exp(y^2/4)\,.
  \end{equation}
We shall refer to time-independent solutions of this equation as expanders.
Eq.\eqref{eq_F} have been very recently studied by Germain and Rupflin \cite{gr} who established interesting results concerning  existence, uniqueness, and stability of expanders. Below we complement these results by a more detailed formal quantitative analysis (which is essential for our purposes).

Expanders satisfy the ordinary differential equation
\begin{equation}\label{exp_eq}
    F''+\left(\frac{d-1}{y}+\frac{y}{2}\right)\,F'-\frac{d-1}{2 y^2} \sin(2F)=0\,
\end{equation}
with the regularity condition $F(y)\sim A y$ near $y=0$,
where $A$ is a free parameter.
In contrast to shrinkers, expanders are globally regular for any $A$. This is due to the strong damping term $\frac{y}{2} F'$ in \eqref{exp_eq} which drives $F'(y)$ rapidly to zero as $y\rightarrow \infty$ and guarantees
that $\lim_{y\rightarrow\infty} F(y)$ exists. Let $B=\lim_{y\rightarrow\infty} F(y)-\pi/2$. It is routine to show that $B$ depends continuously on $A$.
In order to get a more precise asymptotic behaviour, we rewrite Eq.\eqref{exp_eq} in the integral form
\begin{equation}\label{hop}
 F'(y)= \frac{d-1}{2} y^{1-d} e^{-y^2/4} \int_0^y s^{d-3} e^{s^2/4} \sin(2F(s)) ds\,,
\end{equation}
and compute the limit
\begin{equation}\label{hop2}
 \lim_{y\rightarrow\infty} y^3 F'(y)= (d-1) \lim_{y\rightarrow\infty} \frac{\int_0^y s^{d-3} e^{s^2/4} \sin(2F(s)) ds}{2 y^{d-4} e^{y^2/4}}= -(d-1) \sin(2B)\,,
\end{equation}
where the last step follows from l'H\^{o}pital's rule.
Therefore, if $B\neq 0$, we have for large $y$
\begin{equation}\label{asym_F}
    F(y)=\frac{\pi}{2}+B+\frac{d-1}{2 y^2} \sin(2B)+\mathcal{O}(y^{-4})\,.
\end{equation}
 We note in passing that Eq.\eqref{B(A)} below implies that  there is an infinite countable subset of parameter values  for which $B(A)=0$ and
 \begin{equation}\label{exp_special}
 F(y)-\pi/2\sim c\, y^{-d} e^{-y^2/4}\qquad\text{as}\quad y\rightarrow\infty\,.
  \end{equation}
  The variational proof (using a renormalized energy) of existence of such rapidly decaying expanders was recently given in \cite{gr}.
  Since these solutions do not seem to participate in the dynamics of blow-up, we do not pursue them here in more detail.

Next, we derive  asymptotic approximations of the function $B(A)$ for small and large arguments. For small $A$ we substitute $F(y)=A\tilde F(y)$ into Eq.\eqref{exp_eq} and take the limit $A\rightarrow 0$ to obtain the linear equation
\begin{equation}\label{exp_tilde}
    \tilde F''+\left(\frac{d-1}{y}+\frac{y}{2}\right)\,\tilde F'-\frac{d-1}{y^2} \tilde F=0\,
\end{equation}
with the initial condition $\tilde F(y)\sim y$ near $y=0$. Clearly, the solution $\tilde F(y)$ is a positive monotonically increasing function converging to a constant at infinity. The explicit solution is
\begin{equation}
\tilde F(y)=y e^{-\frac{y^2}{4}} M\left(\frac{d+1}{2},\frac{d+2}{2},\frac{y^2}{4}\right)\,,
\end{equation}
where $M(a,b,y)$ is the Kummer confluent hypergeometric function.
Using the asymptotic expansion  $M(a,b,x)\sim\frac{\Gamma(b)}{\Gamma(a)} x^{a-b} e^x$ for large $x$ \cite{nist} we get $\tilde F(\infty)=2\Gamma(\frac{d+2}{2})/\Gamma(\frac{d+1}{2})$, thus for small $A$ we have
\begin{equation}\label{smallA}
    B(A)\simeq -\frac{\pi}{2}+\frac{2\Gamma(\frac{d+2}{2})}{\Gamma(\frac{d+1}{2})}\, A\,.
\end{equation}

For large $A$, repeating the argument leading to Eq.\eqref{asymf}, we get for $1/A\ll y\ll 1$
\begin{equation}\label{asymF}
    F(y)\simeq\frac{\pi}{2} + A^{-\frac{d-2}{2}} \alpha \, y^{-\frac{d-2}{2}} \sin(\omega \ln{y} +\omega \ln{A}+\delta)\,.
\end{equation}
On the other hand, for $y\gg1/A$ we can write
\begin{equation}\label{scal2}
    F(y)\simeq \frac{\pi}{2} +  H(y)\,,
\end{equation}
where $H(y)$ is a solution of the linearized equation
\begin{equation}\label{eq_h}
    H''+\left(\frac{d-1}{y}+\frac{y}{2}\right) H'+\frac{d-1}{y^2} \,H=0
\end{equation}
satisfying $H(\infty)=B$. In contrast to shrinkers, this normalization condition does not determine the solution uniquely since the
 two linearly independent solutions at infinity are
\begin{equation}\label{H12}
    H_1(y)\sim 1 \quad \text{and}\quad H_2(y)\sim y^{-d} \exp(-y^2/4)\,,
\end{equation}
hence
\begin{equation}\label{H}
   H(y)= B H_1(y)+ c \,H_2(y)\,,
\end{equation}
where $c$ is an arbitrary constant.
For $1/A\ll y\ll 1$ the solutions $H_1$ and $H_2$  behave as
\begin{equation}\label{h12asym}
  H_i(y) \simeq  C_i\, y^{-\frac{d-2}{2}} \sin(\omega \ln{y} + \Delta_i)\,,\quad i=1,2.
\end{equation}
Combining Eqs.\eqref{asymF},\eqref{H}, and \eqref{h12asym} we get
the following matching condition
\begin{equation}\label{match_F}
 A^{-\frac{d-2}{2}} \alpha\,  \sin(\omega \ln{y} +\omega \ln{A}+\delta) \simeq B C_1  \sin(\omega \ln{y} + \Delta_1) + c\, C_2\sin(\omega \ln{y} + \Delta_2)\,,
\end{equation}
which yields
\begin{equation}\label{B(A)}
    B(A) \simeq \tilde C A^{-\frac{d-2}{2}}  \sin(\omega \ln{A}+\tilde \delta)\,,
\end{equation}
where $\tilde C$ and $\tilde \delta$ are determined by $\alpha, \delta, C_i,\Delta_i$.

\section{Linear stability of self-similar solutions}
Now, we turn our attention to the linear stability analysis of shrinkers and expanders. The results of this analysis are important in understanding the dynamics of blow-up.
\subsection{Shrinkers}
\noindent Substituting $f(s,y)=f_n(y)+w(s,y)$ into Eq.\eqref{eq_f} and retaining only linear terms in $w$, we get the evolution equation for
linearized perturbations around the shrinker $f_n$
\begin{equation}\label{eq_w}
  w_s =
  \frac{1}{\rho} \left(\rho
  w_y\right)_y -\frac{d-1}{y^2}\,\cos(2f_n)\,w\,,
  \end{equation}
which after separation of variables,
$w(s,y)=e^{-\lambda s} v(y)$, yields the eigenvalue problem
\begin{equation}\label{eq:1}
  \mathcal{A}_n v=\lambda v\,, \qquad \mathcal{A}_n=-\frac{1}{\rho} \partial_y \left(\rho
  \partial_y\right) + \frac{d-1}{y^2} \cos(2f_n)\,.
  \end{equation}
For each $n$ the operator $\mathcal{A}_n$ is self-adjoint in the Hilbert space $X=L_2([0,\infty),\rho\, dy)$. Both endpoints $y=0$ and $y=\infty$ are of the limit-point type with admissible solutions behaving as $v(y)\sim y$ for $y\rightarrow 0$ and $v(y)\sim y^{2\lambda}$ for $y\rightarrow \infty$. Note that for each $n$ there is an eigenvalue $\lambda=-1$ with the associated eigenfunction $v(y)=y f'_n(y)$. The presence of this gauge mode is due to time translation symmetry. To see this observe that if the blow-up time is shifted from $T$ to $T+2\varepsilon$, then
 \begin{equation}\label{T_shr}
    f(y)\rightarrow f\left(\frac{y}{\sqrt{1+2\varepsilon e^s}}\right)=f(y)-\varepsilon e^{s} y f'(y)+\mathcal{O}(\varepsilon^2).
 \end{equation}
 Since $f'_n(y)$ has $(n-1)$ zeroes, it follows from the Sturm oscillation theorem  that for the $n$-th shrinker there are exactly $(n-1)$  eigenvalues below $-1$. We checked numerically (but were unable to prove analytically) that there are no eigenvalues in the interval $-1<\lambda\leq 0$.  Denoting the spectrum by $\{\lambda^{(n)}_k | k=0,1,\dots\}$ we thus have
\begin{equation}\label{spec}
    \lambda^{(n)}_0<\lambda^{(n)}_1<\cdots<\lambda^{(n)}_{n-1}=-1<0<\lambda^{(n)}_n<\cdots
\end{equation}
Concluding, the shrinker $f_n$ has exactly $(n-1)$ unstable modes (the gauge mode with $\lambda=-1$ is not counted as a genuine instability). In particular, the shrinker $f_1$ is linearly stable and therefore it is expected to participate in the generic dynamics of blow-up.
This expectation will be confirmed numerically in section~6. The first few eigenvalues of the operator $\mathcal{A}_n$ for several $n$ in $d=3$  are displayed in Table~II. Note that the columns in this table converge to
to  limiting values, namely for each integer $m$ we have
\begin{align}\label{eq:3}
 \lim_{n\rightarrow\infty}\lambda_{n+m}^{(n)}=\lambda_{m}\,.
\end{align}
Now, we will show that $\lambda_{m}$ are the eigenvalues of the point spectrum of
the operator
\begin{equation}\label{ainf}
\mathcal{A}_{\infty}=-\frac{1}{\rho} \partial_y \left(\rho
  \partial_y\right) - \frac{d-1}{y^2}\,,
\end{equation}
which is obtained from \eqref{eq:1} by taking the (nonuniform) limit $f_n(y)\rightarrow \pi/2$ as $n\rightarrow\infty$. The potential term in \eqref{ainf} is unbounded from below as $y\rightarrow 0$ and $y=0$ is a limit-circle point, so for $\mathcal{A}_{\infty}$  to be self-adjoint,  we have to specify  an additional boundary condition at $y=0$ (which is usually referred to as the self-adjoint extension). This is done as follows.
The solution of the eigenvalue equation $ \mathcal{A}_{\infty} v=\lambda v$
which is admissible at infinity (i.e., behaving as $v(y)\sim y^{2\lambda}$ for $y\rightarrow\infty$) reads
\begin{equation}\label{tricomi}
    v(y)=y^{1-\frac{d}{2}+i\omega} \,U\left(\frac{1}{2}-\frac{d}{4}+\frac{i\omega}{2}
    -\lambda,1+i\omega,\frac{y^2}{4}\right)\,,
\end{equation}
where $U(a,b,z)$ is the Tricomi confluent hypergeometric function.
Using the asymptotic expansion formula for $z\rightarrow 0$ (which is valid for $1\leq\text{Re}(b)<2$) \cite{nist}
\begin{equation}\label{Uasym}
    U(a,b,z)\sim \frac{\Gamma(1-b)}{\Gamma(a-b-1)}+\frac{\Gamma(b-1)}{\Gamma(a)} \, z^{1-b}\,,
\end{equation}
we get from \eqref{tricomi}
\begin{equation}\label{y0asym}
    v(y)\sim y^{1-\frac{d}{2}} \cos\left(\omega\ln{y} +\Phi(\lambda)\right) \qquad \text{as} \,\, y\rightarrow 0\,,
\end{equation}
where
\begin{equation}\label{faza}
    \Phi(\lambda)=\arg\left(\frac{\Gamma(i\omega)}{\Gamma\left(\frac{1}{2}
    -\frac{d}{4}+\frac{i\omega}{2}-\lambda\right)}\right)\,.
\end{equation}
The self-adjoint extension amounts to fixing the phase $\Phi(\lambda)$ modulo $\pi$. A natural choice is to require that the eigenvalue $\lambda=-1$  belongs to the spectrum of $\mathcal{A}_{\infty}$. This leads to the quantization condition
\begin{align}
  \label{eq:7}
  \Phi(\lambda_{m-1})=\Phi(-1)+m\pi,\quad m\in\mathbb{Z}.
\end{align}
As shown in Table \ref{table:A}, solutions of this equation, in fact,  give the limit of the point spectra of the operators $\mathcal{A}_n$ for $n\rightarrow\infty$.

\begin{table}[ht]
  \begin{tabular}{|c|c|c|c|c|c|c|c|c|}\hline
    $\,\,n\,\,$             & $\lambda_{n-4}^{(n)}$ & $\lambda_{n-3}^{(n)}$ &
    $\lambda_{n-2}^{(n)}$ & $\lambda_{n-1}^{(n)}$   & $\lambda_{n}^{(n)}$ &
    $\lambda_{n+1}^{(n)}$    &
    $\lambda_{n+2}^{(n)}$                                                                                             \\\hline
    $1$                     &                         &                         &            & $-1$ & $0.51762$ & $1.63038$ & $2.69684$\\\hline
    $2$                     &                         &                         & $-53.2995$ & $-1$ & $0.48625$ &
    $1.61122$        & $2.68550$                                                                                                       \\\hline
    $3$                     &                         & $-6054.92$              & $-52.4152$ & $-1$ & $0.48271$ & $1.60879$ & $2.68380$\\\hline
    $4$                     & $-699295$               & $-5968.91$              & $-52.3292$ & $-1$ & $0.48237$ & $1.60858$ & $2.68363$\\\hline
    $\vdots$                & $\vdots$                & $\vdots$                &
    $\vdots$                & $-1$                    & $\vdots$                & $\vdots$            & $\vdots$                      \\\hline
    $\infty$                & $-688498$               & $-5959.55$              & $-52.3200$ & $-1$
                            & $0.48234$              & $1.60852$       & $2.68361$                                                    \\\hline
  \end{tabular}
  \caption{The first few eigenvalues of the operator $\mathcal{A}_n$ in $d=3$.
  Numerical solutions of the quantization condition \eqref{eq:7} are listed in the last row.
    \label{table:A}}
\end{table}
\subsection{Expanders}
 \noindent The linear stability analysis of expanders proceeds along the similar lines as above.  Substituting $F(\sigma,y)=F(y)+W(\sigma,y)$ into Eq.\eqref{eq_F} and linearizing we obtain
 the evolution equation for
linearized perturbations around an expander $F(y)$
\begin{equation}\label{eq_W}
  W_{\sigma} =
  \frac{1}{R} \left(R\,
  W_y\right)_y -\frac{d-1}{y^2}\,\cos(2F)\,W\,,
  \end{equation}
which after separation of variables,
$W(\sigma,y)=e^{-\Lambda \sigma} V(y)$, leads to the  eigenvalue problem
\begin{equation}
  \mathcal{B} V=\Lambda V\,, \qquad \mathcal{B}=-\frac{1}{R} \partial_y
   \left(R\, \partial_y
  \right) +\frac{d-1}{y^2} \cos(2F)\,.
  \end{equation}
\addtolength{\topmargin}{-0.5pc}
\addtolength{\textheight}{+2pc}
The operator $\mathcal{B}$ is self-adjoint in the Hilbert space $Y=L_2([0,\infty),R\,dy)$. Both endpoints are of the limit-point type with admissible solutions $V(y)\sim y$ for $y\rightarrow 0$ and $V(y)\sim y^{2\lambda-d} e^{-y^2/4}$ for
$y\rightarrow\infty$. The gauge mode due to time translation symmetry $V(y)=y F'(y)$ has the (formal) eigenvalue $\Lambda=1$, because if $T\rightarrow T+2\varepsilon$, then
 \begin{equation}\label{T_exp}
    F(y)\rightarrow F\left(\frac{y}{\sqrt{1-2\varepsilon e^{-\sigma}}}\right)=F(y)+\varepsilon e^{-\sigma} y F'(y)+\mathcal{O}(\varepsilon^2).
 \end{equation}
 The gauge mode is not an eigenfunction (because it does not belong to $Y$), nevertheless the Sturm oscillation theorem still applies and implies that an
expander with $n$ zeros of $F'(y)$ has exactly $n$ eigenvalues below $+1$ (this was proved independently in \cite{gr}).
In particular, monotonic expanders are linearly stable. Although this fact will be sufficient for the analysis of continuation beyond the generic blow-up, we wish to point out that using the "turning-point" method \cite{s} one can determine sharp stability intervals for expanders. This is done as follows.
Let $F_A(y)$  denote the expander starting with $F'(0)=A$. Differentiating Eq.\eqref{exp_eq} with respect to $A$ we find that
$\frac{\partial F_A(y)}{\partial A}$ is the zero mode of the operator $\mathcal{B}$.
In general, $\frac{\partial F_A(y)}{\partial A}\sim B'(A)\neq 0$ for  $y\rightarrow \infty$, however it follows from Eq.\eqref{B(A)} that there is an increasing sequence of numbers $A_k$ ($k\in \mathbb{}{N}$) for which $B'(A_k)=0$ and then, by Eq.\eqref{exp_special}
\begin{equation}\label{}
    \frac{\partial F_A(y)}{\partial A} \sim c'(A)\, y^{-d} \exp(-y^2/4)\,,
\end{equation}
 hence the zero mode is a genuine eigenfunction. By \cite{s} this implies that $A_k$
 are turning points at which the expander $F_A$  picks a new unstable mode. More precisely,
 the expander $F_A$ with $A\in (A_{k-1},A_k)$ has exactly $(k-1)$ instabilities (here $A_0=0$ by definition).

 \section{Continuation beyond blow-up}
 Suppose that the solution of Eq.\eqref{heat_flat} develops a type I singularity at time $T$. Then, as we showed above, the profile of blow-up is given by one of the shrinkers
\begin{equation}\label{blow1}
    \lim_{t\nearrow T} u(t,r \sqrt{T-t})=f_n(r)\,,\qquad n\in \mathbb{N}\,.
\end{equation}
In order to continue the solution beyond blow-up, for times $t>T$  we glue an expander which matches the shrinker $f_n$ at time $T$, that is we require that
\begin{equation}\label{glue1}
  \lim_{t\searrow T} u(t,r\sqrt{t-T})=F_A(r)\,,\qquad  \hspace{1.3cm} F_A(\infty)=f_n(\infty)\Longleftrightarrow B(A)=b_n\,,
\end{equation}
or
\begin{equation}\label{glue2}
 \lim_{t\searrow T} u(t,r \sqrt{t-T})=\pi-F_A(r)\,,\qquad
\pi-F_A(\infty)=f_n(\infty)\Longleftrightarrow B(A)=-b_n\,,
\end{equation}
Note that in the case \eqref{glue1} the solution  stays continuous across blow-up (hence the degree does not change), while in the case \eqref{glue2} the solution jumps at the origin from $u(t,0)=0$ for $t<T$ to $u(t,0)=\pi$ for $t>T$ (hence the degree changes by one).
 In both cases we obtain a global weak solution which is smooth except for the time $T$.

 Let $N(n)$ denote the number of roots of the equation $|B(A)|=|b_n|$.
 It follows from the large $A$ formula   for expanders \eqref{B(A)} and large $n$ formula for shrinkers \eqref{scala} that $N(n)$ increases indefinitely with $n$. More precisely, we find numerically (see Fig.~3) that for $n\geq 2$
 \begin{equation}
N(n)=
\begin{cases} 2n-3 & \text{for $d=3,4$,}
\\
2n-1 & \text{for $d=5,6$.}
\end{cases}
\end{equation}
Since the shrinker $f_n$ has $(n-1)$ instabilities, all $n\geq 2$ blow-ups are non-generic phenomena of codimension $(n-1)$.
 It follows from the stability analysis of expanders that only one continuation is stable, namely that with $A^*_n= \min\{A:|B(A)|=|b_n|\}$. For this stable continuation the degree changes by one if $n$ is odd and does not change if $n$ is even.

\begin{figure}[ht]
  \centering
  \includegraphics[width=.48\linewidth]{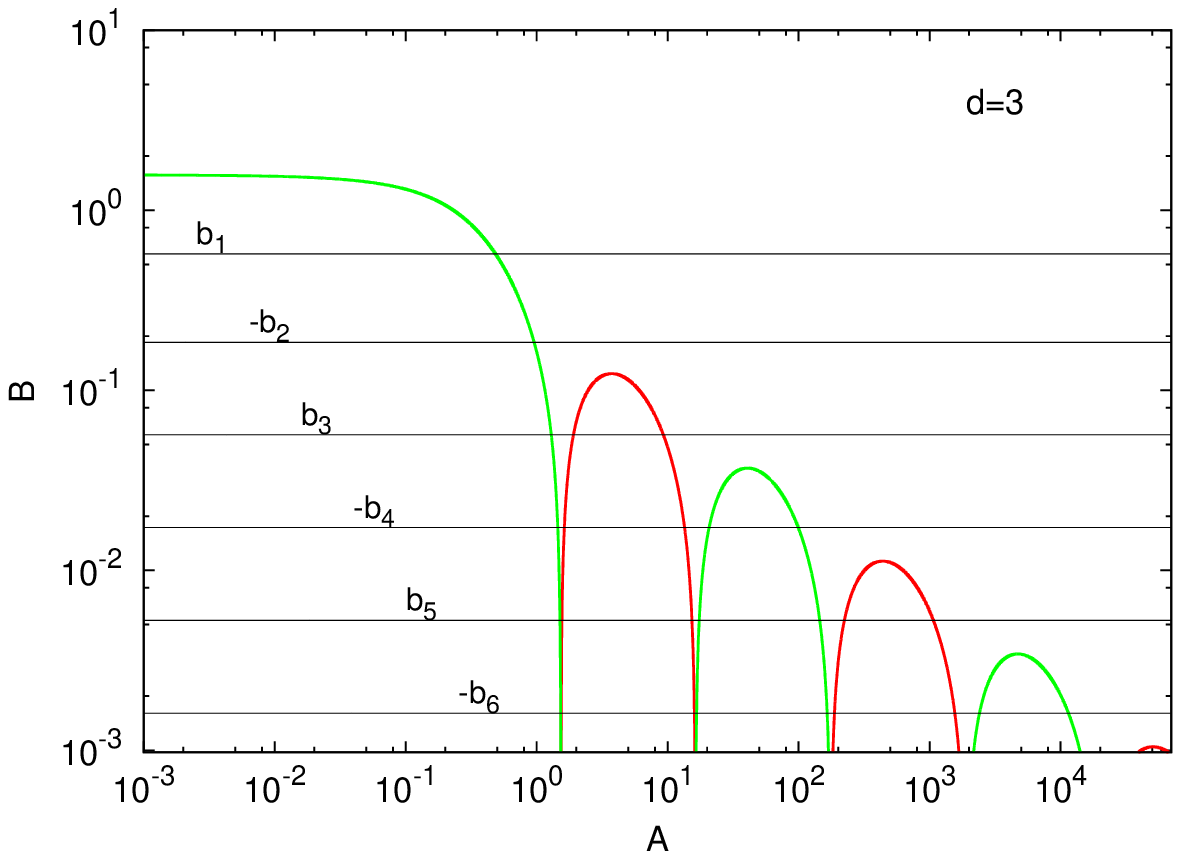}
    \includegraphics[width=.48\linewidth]{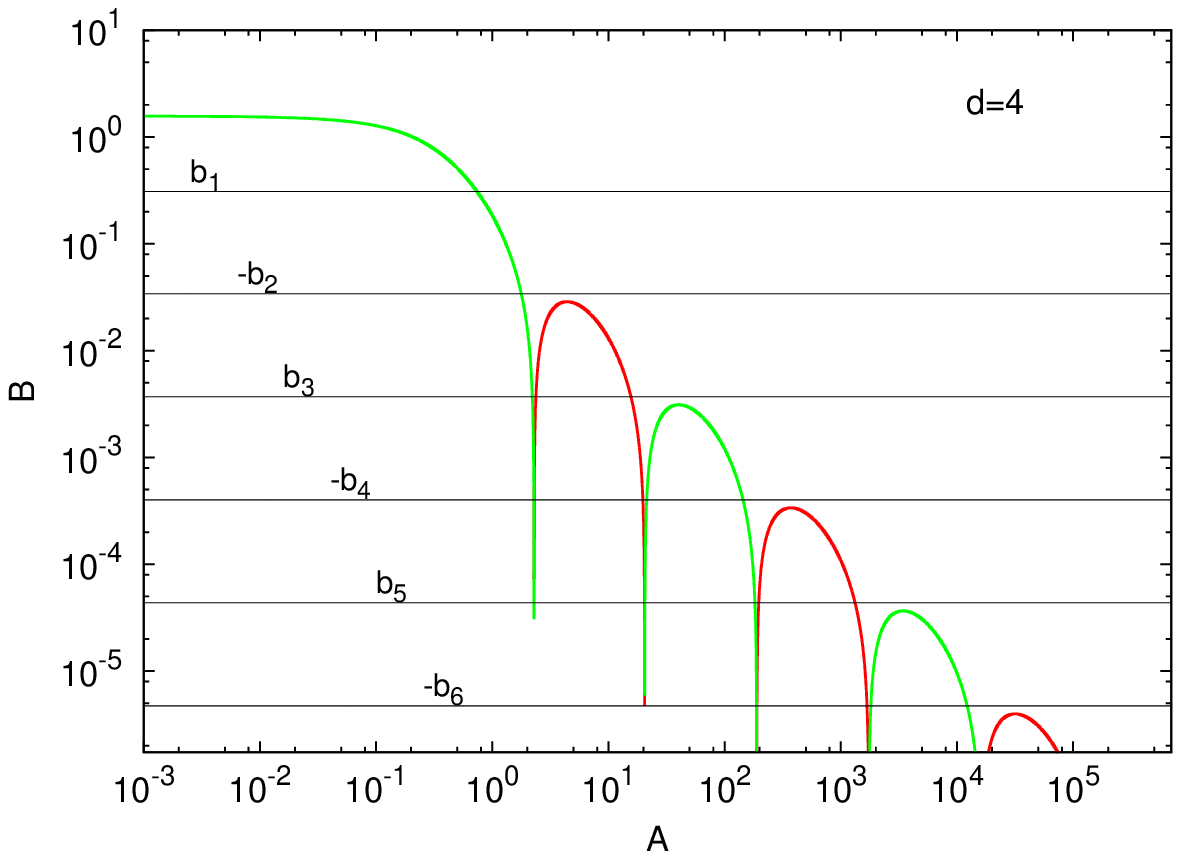}
      \includegraphics[width=.48\linewidth]{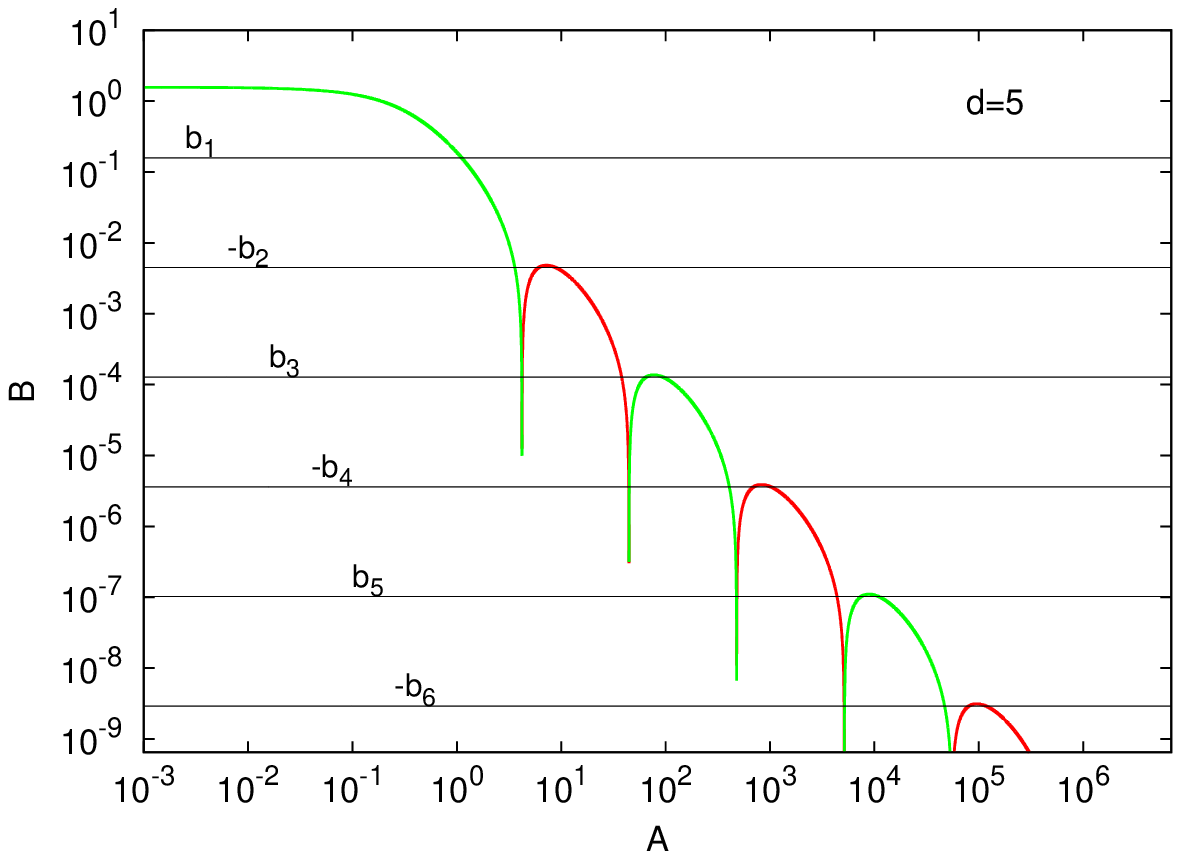}
        \includegraphics[width=.48\linewidth]{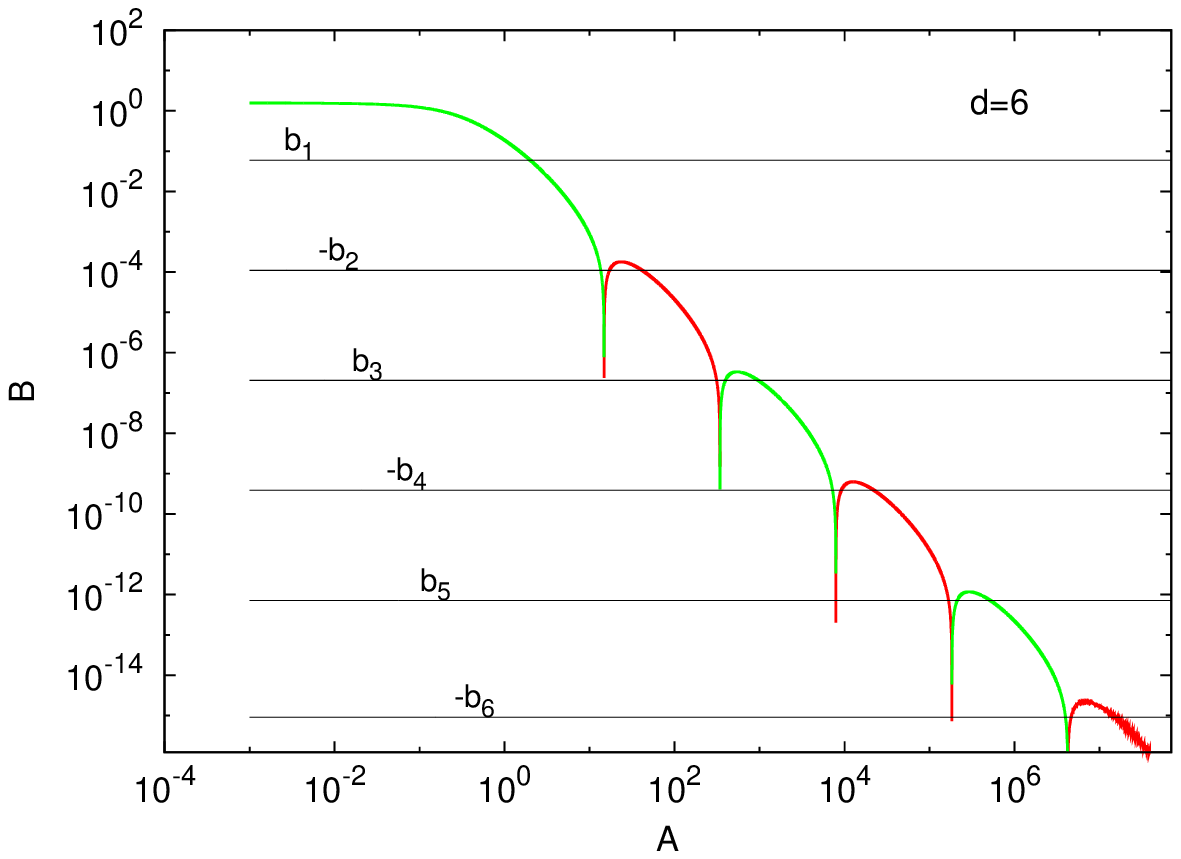}
  \caption{Plots of $|B(A)|$ in the log-log scale. The intersections with horizontal lines $|B|=|b_n|$ determine the number of continuations beyond blowup.}
  \label{fig:roots}
\end{figure}

 Hereafter, we focus on the most important and interesting case $n=1$ corresponding to the generic blow-up governed by the linearly stable shrinker $f_1$.
 In this case the equation $B(A)=b_1$ has no roots, while the equation $B(A)=-b_1$ has exactly one root (note that the existence of this root is guaranteed by the small $A$ formula \eqref{smallA} and the continuity of the function $B(A)$), hence the continuation beyond blow-up is unique, stable, and changing degree.
As this is our main result, let us phrase it in the form of a conjecture:
\vskip 0.2cm
\noindent  \emph{\bf{Conjecture 1}.} Let $3\leq d\leq 6$. Suppose that $u(t,r)$ is
 a generic solution of Eq.\eqref{heat_flat} which develops a singularity at $r=0$ in a finite time $T$. Then, for sufficiently small $r$ there holds
\begin{equation}
u(t,r)\sim
\begin{cases} f_1\left(\dfrac{r}{\sqrt{T-t}}\right) & \text{for $T-r^2<t<T$,}
\\
\dfrac{\pi}{2}+b_1 & \text{for $t=T$,}
\\
\pi-F_{A^*_1}\left(\dfrac{r}{\sqrt{t-T}}\right) & \text{for $T<t<T+r^2$.}
\end{cases}
\end{equation}
where $A^*_1$ is the (unique) root of equation $B(A)=-b_1$. This is illustrated in  Fig.~4.
\vskip 0.2cm
In the next section we present numerical evidence supporting this conjecture.
\begin{figure}[h]
  \centering
  \includegraphics[width=.55\linewidth]{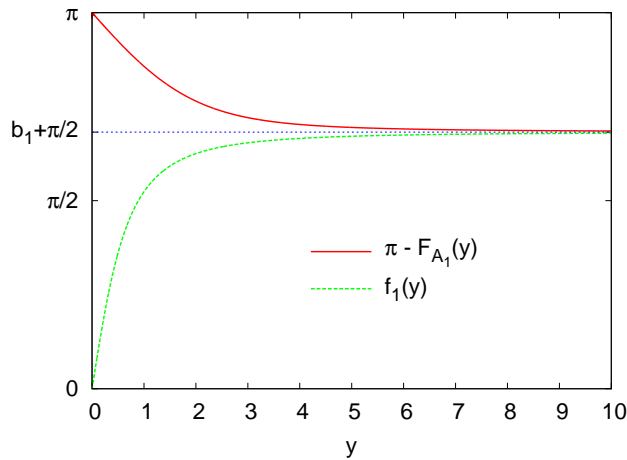}
  \caption{Gluing together the stable shrinker and expander in $d=3$. Here $b_1=0.573141$ and $A^*_1=0.483668$. }
  \label{fig:glue}
\end{figure}
\vspace{-1cm}
\section{Numerical evidence}
In this section we verify the above heuristic
predictions by numerical simulations. In order to keep track of the structure of the singularity developing on a vanishingly small scale, it is necessary to use an adaptive method which refines the spatio-temporal grid near the singularity.
 Our numerical method is based on
the moving mesh method combined with  the Sundman transformation, as described
in \cite{bw}, with some minor modifications and improvements specific to the problem at hand. This method is particularly efficient in computations of self-similar singularities.
To implement the adaptivity in time we introduce a new computational time variable $\tau$ defined by
\begin{equation}\label{sund}
    \frac{dt}{d\tau}=g(u)\,,\quad g(u)=\left\lvert\frac{u_r}{u_{rt}}\right\rvert_{r=0}\,.
\end{equation}
Under this rescaling (called the Sundman transformation) the fixed time steps in $\tau$
correspond to $\Delta t_i\approx (T-t_i) \Delta\tau$ as $t\nearrow T$. In this way, the time scale of the developing singularity is identified automatically even though the blow-up time $T$ is unknown beforehand.
To implement the adaptivity in space we introduce a new computational spatial variable $\xi\in[0,1]$ and define a mesh function $r(\xi,t)$ which places the moving mesh points at $r_i(t)=r(i\Delta\xi,t)$. The function  $r(\xi,t)$, whose role is to cluster the mesh points near the singularity, is determined by an auxiliary moving mesh partial differential equation (MMPDE), which is solved simultaneously with the original PDE. We use the so called MMPDE6 \cite{hrr}
\begin{equation}
  \label{eq:9}
 \varepsilon r_{t\xi\xi}= -(M r_\xi)_\xi\,,
 \end{equation}
 with the mesh density
  $M=|u_r|+\sqrt{|u_{rr}|}$ and the time-dependent relaxation parameter $\varepsilon(t)=100\sqrt{g(t)}+0.05$ (this $\varepsilon(t)$, found empirically, results in a better performance than the
customarily used constant value).

The harmonic map heat equation $u_t=N(u)$, where $N(u)$ is the right hand side of Eq.\eqref{heat} or Eq.\eqref{heat_flat}, is  now be rewritten as the system
\begin{eqnarray}
    t_\tau&=&g(u)\,,\\
    u_\tau+r_\tau u_r&=&g(u) N(u)\\
    \varepsilon r_{\tau\xi\xi}&=&-g(u)(M r_\xi)_\xi\,.
\end{eqnarray}
These equations are discretized using a 5-point finite difference scheme and integrated
via the Embedded Runge-Kutta-Fehlberg (RKF45) method.

The numerical results are presented below for $d=3$ as an illustration; the behaviour of solutions is qualitatively the same in all dimensions $3\leq d\leq 6$. Since the dynamics of blow-up does not depend on the curvature of the domain, we first show simulations for Eq.\eqref{heat_flat}, and only at the end we show simulations of multiple blow-ups for the spherical domain Eq.\eqref{heat}.
%

We begin by demonstrating the convergence to the stable shrinker. Fig.~5 depicts snapshots from a typical evolution ending in a singularity. As the blow-up is approached, the solution is seen to converge to the profile of the stable shrinker $f_1$.
\begin{figure}[h]
  \centering
  \includegraphics[width=.98\linewidth]{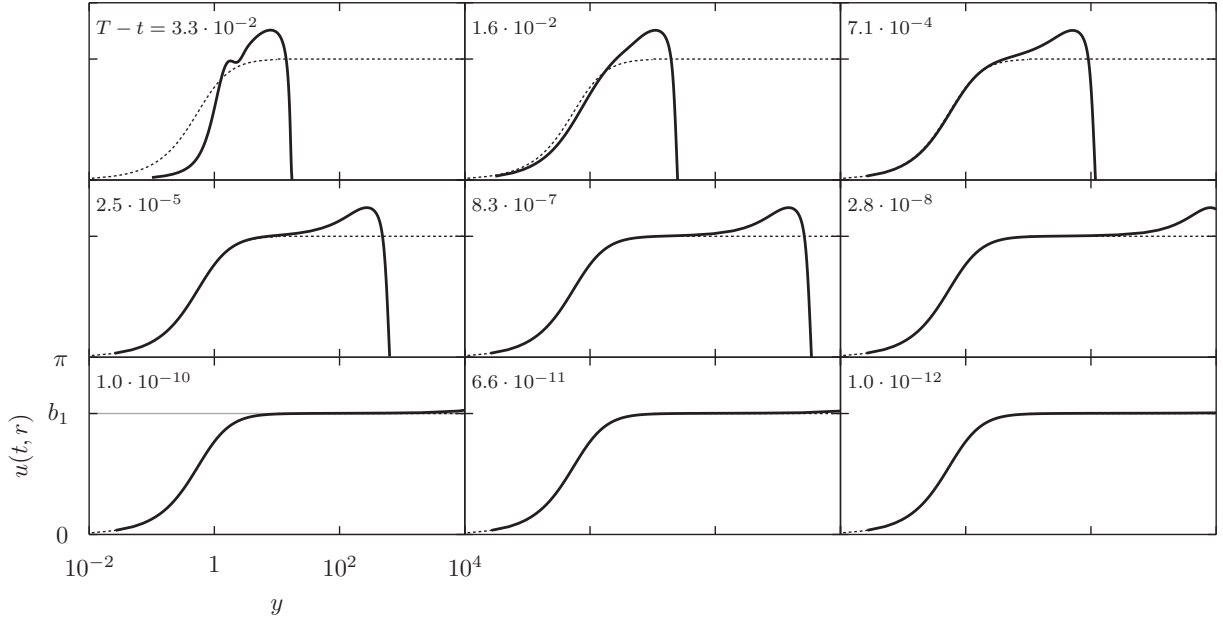}\\
  \caption{Convergence of the numerical solution (solid line) to the stable shrinker $f_1(y)$
    (dotted line).}
  \label{fig:conv_to_f1}
  \end{figure}

According to the linearized stability analysis the deviation of the solution from the stable shrinker is expected to have the following form near $r=0$ for $t\nearrow T$ (to avoid notational clutter, hereafter  we drop the superscript $(1)$ on the eigenvalues and the eigenfunctions)
\begin{equation}\label{blow_expand}
    u(t,r)-f_1(y)\simeq \sum_{k=1}^{\infty} c_k (T-t)^{\lambda_k}\, v_k(y) = c_1
    (T-t)^{\lambda_1}\, v_1(y)+\mathcal{O}\left((T-t)^{\lambda_2}\right)\,,
\end{equation}
where $y=r/\sqrt{T-t}$.
To verify this prediction we proceed as follows. Differentiating \eqref{blow_expand} twice and using the normalization $v_1'(0)=1$ we obtain
\begin{equation}
  \label{eq:2}
  \partial_t\left[(T-t)^{1/2}\, u_r\right] \big\rvert_{r=0}
  =-c_1\lambda_1(T-t)^{\lambda_1}+\mathcal{O}\left((T-t)^{\lambda_2}\right)\,.
\end{equation}
Fitting the right-hand side of this equation to the numerically computed left-hand side,
we get the coefficient $c_1$ and the eigenvalue $\lambda_1$ (see the left panel of Fig.~6). The fit gives $\lambda_1=0.519$, in good accord with the linearized stability analysis (see Table~II). Next, in the right panel of Fig.~6 we show that near the blow-up time the left- and the right-hand sides of the expression \eqref{blow_expand}
(computed completely independently) do indeed agree.
\addtolength{\topmargin}{-1pc}
\addtolength{\textheight}{1.6pc}
\begin{figure}[h]
  \centering
  \includegraphics[width=\linewidth]{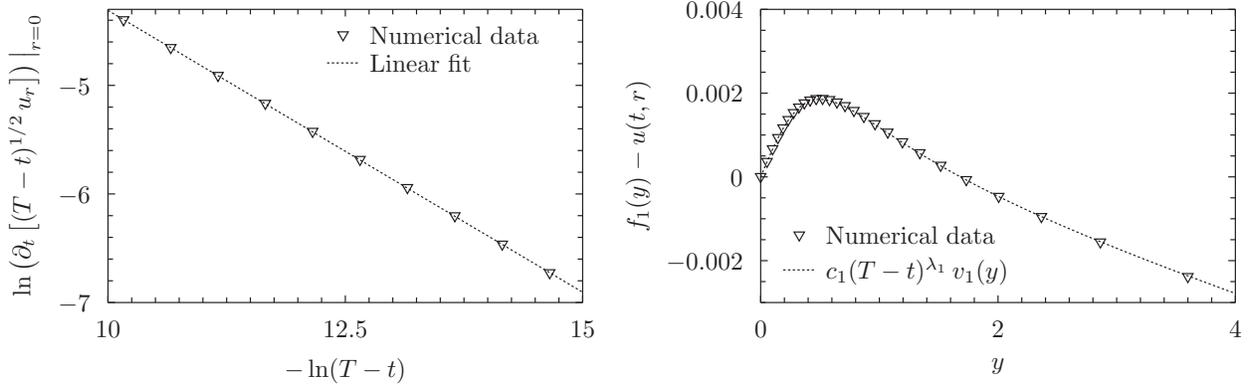}
  \caption{Left: The log-log plot of the left-hand side of expression \eqref{eq:2}. The linear fit  gives $\lambda_1=0.519$. Right: We plot the deviation of the numerical solution from the stable shrinker at $T-t=3.19\cdot 10^{-6}$ and superimpose the first stable eigenmode $c_1 (T-t)^{\lambda_1} v_1(y)$, obtained by solving the eigenvalue equation \eqref{eq:1}, with the coefficient $c_1$ taken from the fit in the left panel.}
  \label{fig:first_mode}
\end{figure}

Next, we describe the continuation beyond blow-up. In order to pass through the singularity we need to modify the numerical code.
First, according to \eqref{glue2} we expect that at the blow-up time the solution
is discontinuous at $r=0$. This behaviour  is not compatible
with the boundary condition $u(t,0)=0$ implemented in our code. To go around this difficulty, we simply
 rewrite Eq.\eqref{heat_flat} in terms of
$z(t,r)=ru(t,r)$ and impose the boundary condition $z(t,0)=0$ (which is compatible with the jump). In the case of Eq.\eqref{heat}  we use a similar trick introducing
$Z(t,\theta)=\sin(\theta)\, U(t,\theta)$ as an independent variable.
Second, at some late stage of blow-up
(say, $T-t=10^{-10}$) we must switch off the Sundman transformation since otherwise
the time step would keep decreasing down to the machine precision, effectively freezing the simulation and preventing  it to cross the time of blow-up. To this end, we replace
$g(u)$ in \eqref{sund} by $G(u)=g(u)+\Delta$ where $\Delta\approx 10^{-10}$ serves as a small scale cut-off.
When $g(u)\ll\Delta$, the solver loses its ability to adapt the
time step appropriately and very quickly steps over the
blow-up time. A moment afterwards, when $g(u)$ exceeds $\Delta$ again, the
Sundman transformation is turned back on and keeps tracking of the, now growing,
time-scale of the expander. During a short time interval $T-10^{-10}\lesssim t\lesssim T+10^{-10}$ when the time adaptation procedure is suspended, the spatio-temporal scales are unresolved and the numerical solution is inaccurate (the third and the fourth snapshot in Fig.~7).

\begin{figure}[h]
  \centering
  \includegraphics[width=.98\linewidth]{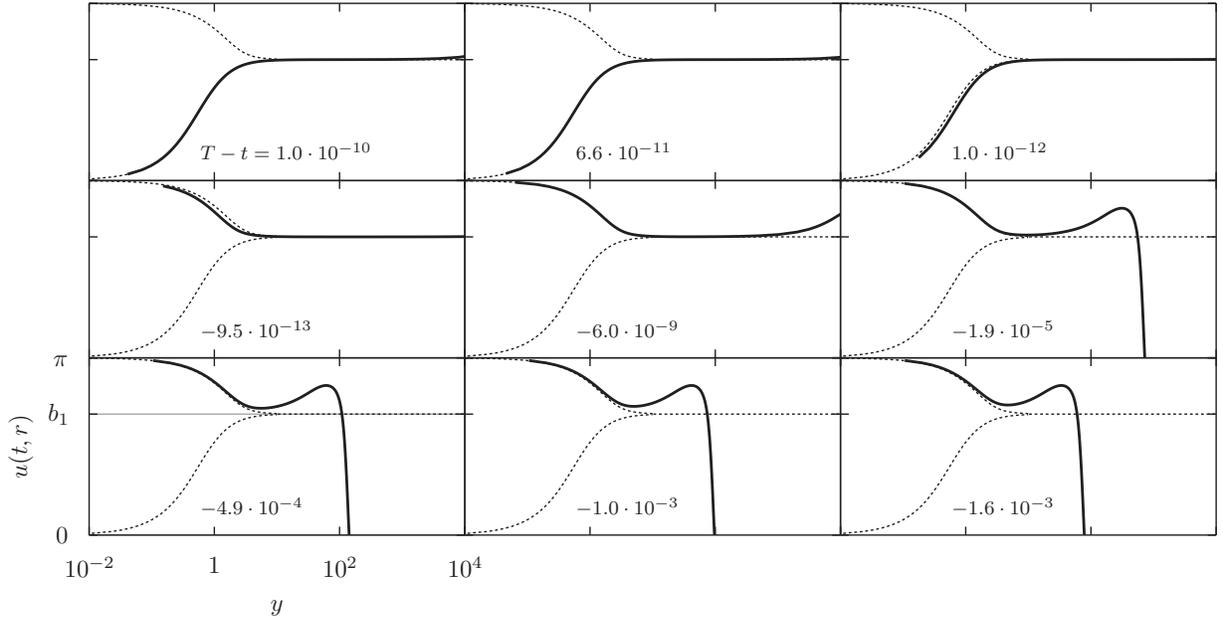}
  \caption{The same evolution as in Fig.~\ref{fig:conv_to_f1} but
    using the modified numerical method which allows the solution to pass
    through the singularity. An additional dotted line shows the
    expander $\pi-F_{A^*_1}(y)$. Notice that for $|T-t|\lesssim 10^{-11}$ the spatio-temporal
    resolution is lost and the numerical solution slightly deviates from the shrinker $f_1(y)$ (the third snapshot) and the expander $\pi-F_{A^*_1}(y)$ (the fourth snapshot). When the resolution is regained, the solution
    converges to the expander but later it moves away from it due to the interference with the far-field structure.}
  \label{fig:single_blow-up}
\end{figure}

 \addtolength{\topmargin}{1pc}
\addtolength{\textheight}{-2pc}
Applying this method, we continue the evolution shown in Fig.~\ref{fig:conv_to_f1}
 past the singularity. In accord with Conjecture~1, almost immediately after the blow-up the numerical solution takes the form of the expander $\pi-F_{A^*_1}(y)$ (see Fig.~7). As written above,  numerical evolution  through a singularity necessarily involves an interval of uncontrolled behaviour due to the inevitable loss of resolution near the instant of blow-up. For this reason the simulation has limited reliability and taken alone would not provide ample evidence for the conjectured behaviour. It is the excellent consistency  between numerics
 and the analytic insight, based on the understanding of self-similar solutions and their linear perturbations, which makes us feel confident that our conjecture is true.

Finally, let us consider the heat flow for harmonic maps between spheres $U:S^d\rightarrow S^d$. As emphasized above, the curvature of the domain manifold
 is  irrelevant in the formation of point singularities, hence all the above results concerning the asymptotic dynamics of blow-up (in particular Conjecture~1) remain valid in the case of a spherical domain. What makes the spherical domain interesting is a pattern of multiple blow-ups for high-degree initial maps. This is illustrated in Fig.~8  showing three consecutive blow-ups at the north pole, south pole, and again the north pole (animated simulations can be found at \cite{pb}). At each blow-up the degree of the map changes by one and eventually the solution comes to rest at the zero energy constant map. Note that, in view of
 the monotonicity formula \eqref{dEdt} and Struwe's theorem (asserting that for harmonic maps between  compact manifolds  the heat flow starting from an initial map with nonzero degree and sufficiently small energy must blow up in finite time), Conjecture~1 implies that the solution starting from an initial map of degree $k$ must blow-up  at least $k$ times (note that the degree of the map need not decrease monotonically).

\begin{figure}[h]
  \centering
  \includegraphics[width=.67\linewidth]{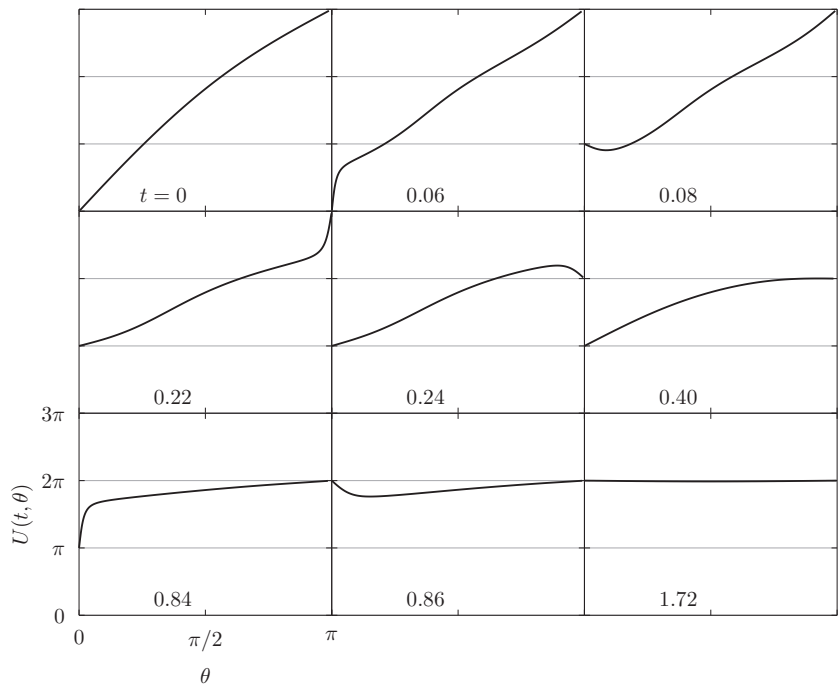}
  \caption{The solution of Eq.\eqref{heat} starting from the initial map $U_0(\theta)=\sin{\theta}+3\theta$ of degree $3$.
  After three blow-ups the map becomes topologically trivial and settles down to the constant map.}
  \label{fig:winding2}
\end{figure}
\vspace{-1cm}
 \addtolength{\topmargin}{-0.5pc}
\addtolength{\textheight}{1pc}
\section{Final remarks}
As mentioned in the introduction, the global weak peaking solutions (having the form of a shrinker and an expander glued together at infinity) exist for many supercritical heat flow equations, so it is natural to ask if these equations, similarly to the harmonic map flow, enjoy the uniqueness of continuation beyond the generic blow-up. We are currently investigating this question in the following models:
\begin{itemize}

\item \textbf{$\ell$-equivariant harmonic map flow:}  The corotational ansatz $(r,\phi)\rightarrow (u(r),\phi)$ is the special ($\ell=1$) case of a more general $\ell$-equivariant ansatz $(r,\phi)\rightarrow (u(r),\chi_{\ell} (\phi))$, where $\chi_{\ell}:S^{d-1}\rightarrow S^{d-1}$ is an eigenmap with constant energy density $k=\ell(\ell+d-2)/2$. For $\ell$-equivariant maps
 Eq.\eqref{heat_flat} changes to
   \begin{equation}\label{heat_ell}
  u_t=  \frac{1}{r^{d-1}} \left(r^{d-1} u_r \right)_r- \frac{k}{r^2} \sin(2u)\,.
\end{equation}
All the qualitative results concerning existence of shrinkers and expanders and their linear stability obtained above for $\ell=1$ trivially carry over to $\ell>1$ provided that $3\leq d <2\ell+2\sqrt{\ell}+2$, however the quantitative characteristics of self-similar solutions (in particular those which imply the uniqueness of  gluing an expander to the stable shrinker) remain to be checked.

\item \textbf{Yang-Mills heat flow:}  It is well-known that there are close parallels between the harmonic map and the Yang-Mills heat flows \cite{g2}. For the spherically symmetric magnetic Yang-Mills potential $h(t,r)$ in $d\geq 3$ dimensions  the analogue of Eq.\eqref{heat_flat} reads
     \begin{equation}\label{ym}
  h_t=  \frac{1}{r^{d-3}} \left(r^{d-3} h_r \right)_r- \frac{d-2}{r^2}\, h (h-1) (h-2)\,.
\end{equation}
Using a similar shooting technique as in \cite{f} one can easily show that for $5\leq d\leq 9$ there are infinitely many shrinkers $h(t,r)=\phi_n(y)$. One novel feature, in comparison with the harmonic map flow, is that the first (stable) shrinker is known explicitly \cite{w}:
\begin{equation}
\phi_1(y)=\frac{y^2}{b+ay^2},\quad
b=\frac{1}{2}(6d-12-(d+2)\sqrt{2d-4})\,,\quad a=\frac{\sqrt{d-2}}{2\sqrt{2}}\,.
\end{equation}
This may be helpful in proving the Yang-Mills analogue of Conjecture~1.

\item \textbf{Semilinear heat equation:} The equation
\begin{equation}\label{power}
u_t=\Delta u + |u|^{p-1} u
\end{equation}
for $d\geq 3$ and supercritical powers
\begin{equation}
\frac{d+2}{d-2}<p<p^*:=
\begin{cases} \infty & \text{for $3\leq d \leq 10$},
\\
 1+\frac{6}{d-10} & \text{for $d\geq 11$,}
\end{cases}
 \end{equation}
has self-similar solutions (shrinkers and expanders) \cite{l,hw,hy} which give rise to global peaking solutions similar to the ones described in section~5, however all these solutions are unstable \cite{gv} (cf. also \cite{mm,fmp}).
 It seems interesting to see if a kind of analogue of Conjecture~1 holds for (codimension-one) threshold solutions.
\end{itemize}

In this paper we restricted our analysis to dimensions $3\leq d\leq 6$. We wish to emphasize that this is not a technical restriction. For $d\geq 7$ the shrinkers disappear and consequently the blow-up changes character from type I to type II \cite{bbv}.

\vskip 0.2cm \noindent \textbf{Acknowledgments:} The second author acknowledges discussions with Marek Fila, Pierre Germain, and Michael Struwe. Special thanks are due to Juan Vel\'{a}zquez for very helpful remarks and suggestions.
The authors are grateful to the Erwin Schr\"odinger Institute in Vienna, where part of this work was done in February 2010 during the program "Quantitative Studies
of Nonlinear Wave Phenomena". This work was supported by the Foundation for Polish Science under the MPD Programme "Geometry and Topology in Physical Models"
co-financed by the EU European Regional Development Fund.


\begin{thebibliography}{10}

\bibitem{bc} P. Bizo\'n and T. Chmaj,  \emph{Harmonic maps between spheres},  Proc. Roy. Soc. London Ser. A  453, 403--415 (1997)

 \bibitem{cw} K. Corlette and R.M. Wald, \emph{Morse theory and infinite families of harmonic maps between spheres},  Comm. Math. Phys. 215, 591--608 (2001)


\bibitem{f} H. Fan,  \emph{Existence of the self-similar solutions in the heat flow of harmonic maps},  Sci. China Ser. A  42, 113-132  (1999)

\bibitem{i} T. Ilmanen, \emph{Lectures on mean curvature flow and related equations}, Lecture Notes,
ICTP, Trieste, 1995, http://www.math.ethz.ch/˜ilmanen/papers/pub.html

\bibitem{lt} A.A. Lacey and D.E. Tzanetis, \emph{Global, unbounded solutions to a parabolic equation}, J. Differential Equations 101, 80-102 (1993)

\bibitem{gv} V.A. Galaktionov and J.L. V\'{a}zquez, \emph{Continuation of blowup solutions of nonlinear heat equations in several space dimensions}, Comm. Pure Appl. Math. 50, 1-67 (1997)

\bibitem{h} G. Huisken, \emph{Asymptotic behavior for singularities of the mean curvature flow}, J. Differential Geom. 31, 285-299 (1990)

\bibitem{aci} S.B. Angenent, D. Chopp, and T. Ilmanen, \emph{A computed example of nonuniqueness
of mean curvature flow in $\mathbb{R}^3$}, Comm. Partial Differential Equations 20, 1937-1958 (1995)

\bibitem{aiv} S.B. Angenent, T. Ilmanen, and J.J.L. Vel\'{a}zquez, \emph{Fattening from smooth initial data
in mean curvature flow}, preprint

\bibitem{gas} A. Gastel, \emph{Nonuniqueness for the Yang-Mills heat flow}, J. Differential Equations 187, 391–411 (2003)

\bibitem{fik} M. Feldman, T. Ilmanen, and D. Knopf, \emph{Rotationally symmetric shrinking and expanding gradient K\"{a}hler-Ricci solitons},  J. Differential Geom.  65, 169-209  (2003).

\bibitem{gal2} V.A. Galaktionov, \emph{Incomplete self-similar blow-up in a semilinear fourth-order reaction-diffusion equation}, 	arXiv:0902.1090

\bibitem{s1} M. Struwe, \emph{Geometric evolution problems},  in: Nonlinear PDE in differential geometry (Park City, UT, 1992), 257–339, IAS/Park City Math. Ser., 2, Amer. Math. Soc., Providence, RI.

\bibitem{gr} P. Germain and M. Rupflin, \emph{Self-similar expanders of the harmonic map flow},  arXiv:1010.6259 [math.AP]

\bibitem{nist} NIST Handbook of Mathematical Functions, edited by F.W.J. Olver, D.W. Lozier, R.F. Boisvert, and C.W. Clark, Cambridge University Press, Cambridge, 2010.

\bibitem{s} R. Sorkin, \emph{A Criterion for the Onset of Instability at a Turning Point}, Astrophysical  J. 249, 254-257 (1981)

\bibitem{bw} C.J. Budd and J. F. Williams, \emph{How to adaptively resolve evolutionary singularities in differential equations with symmetry},
    Journal of Engineering Mathematics 66, 217-236 (2010)

\bibitem{hrr} W. Huang, Y. Ren, and R.D. Russell, \emph{Moving mesh partial differential equations (MMPDES) based on the equidistribution principle}, SIAM J. Numer. Anal. 31, 709-730 (1994)

\bibitem{pb} \url{http://th.if.uj.edu.pl/~biernat/movies}

\bibitem{g2} A. Gastel, \emph{Singularities of first kind in the harmonic map
and Yang-Mills heat flows}, Math. Z. 242, 47–62 (2002)

\bibitem{w} B. Weinkove, \emph{Singularity formation in the Yang-Mills flow}, Calc. Var. Partial Differential Equations 19, 211-220 (2004)

\bibitem{l} L.A. Lepin, \emph{Self-similar solutions of a semilinear heat equation,} Mat. Model. 2, 63–74 (1990)

\bibitem{hw}  A. Haraux and F. B. Weissler, \emph{Nonuniqueness for a semilinear initial value problem,}
Indiana Univ. Math. J. 31, 167-189 (1982)

\bibitem{hy} M. Hirose and E. Yanagida, \emph{Global Structure of Self-Similar Solutions in a Semilinear Parabolic Equation},
J. Math. Analysis and Applications 244, 348-368 (2000)

\bibitem{mm} H. Matano and F. Merle, \emph{On nonexistence of type II blowup for a supercritical nonlinear heat equation}, Comm. Pure and Applied Math. LVII, 1494-1541 (2004)

\bibitem{fmp} M. Fila, H. Matano, and P. Pol\'{a}\v{c}ik, \emph{Immediate regularization after blow-up}, SIAM J. Math. Anal. 32, 752-776 (2005)

\bibitem{bbv} P. Biernat, P. Bizo\'n,  and J.L.L. Vel\'{a}zquez, in preparation


\end{thebibliography}
\end{document}